\documentclass[11pt]{article}

\usepackage{amssymb}
\usepackage{amsmath}
\usepackage{theorem}
\usepackage{epsfig}
\usepackage{verbatim}
\usepackage{graphicx}

\newtheorem{thm}{Theorem}[section]

\newtheorem{prop}[thm]{Proposition}
\newtheorem{defi}[thm]{Definition}

\newcommand{\R}{\Bbb{R}}

\newcommand{\T}{\mathbb{T}}
\newcommand{\D}{\displaystyle}

\newcommand{\grad}{\nabla}

\newcommand{\dt}{\frac{d}{dt}}

\newcommand{\gradp}{\grad^{\bot}}

\newcommand{\dxu}{\partial_{x_1}}
\newcommand{\dxd}{\partial_{x_2}}

\newcommand{\dpt}{\partial_t}

\newcommand{\la}{\lambda}

\newcommand{\al}{\alpha}
\newcommand{\ep}{\varepsilon}
\newcommand{\G}{\Gamma}
\newcommand{\g}{\gamma}
\newcommand{\dg}{\partial_{\gamma}}
\newcommand{\e}{\eta}
\newcommand{\de}{\partial_{\eta}}

\begin{document}

\author{Francisco Gancedo}
\title{Existence for the $\al$-patch model and the\\ QG sharp front in Sobolev spaces}

\date{}

\maketitle

\begin{abstract}
We consider a family of contour dynamics equations depending on a parameter $\al$ with
$0<\alpha\leq 1$. The vortex patch problem of the 2-D Euler equation is obtained taking
$\alpha\rightarrow 0$, and the case $\alpha=1$ corresponds to a sharp front of the QG equation. We
prove local-in-time existence for the family of equations in Sobolev spaces.

\end{abstract}

\maketitle


\section{Introduction}

The 2-D QG equation provides particular solutions of the evolution of the temperature from a
general quasi-geostrophic system for atmospheric and oceanic flows. This equation is derived
considering small Rossby and Ekman numbers and constant potential vorticity (see \cite{P} for more
details). It reads

\begin{align}
\begin{split}\label{QG}
\theta_t(x,t)&+u(x,t)\cdot\grad\theta(x,t)=0,\quad x\in\R^2,\\
\theta(x,0)&=\theta_0(x).
\end{split}
\end{align}
Here $\theta$ is the temperature of the fluid, the incompressible velocity $u$ is expressed by
means of the stream function as follows
\begin{equation*}
u=\gradp\psi=(-\dxd\psi,\dxu\psi),
\end{equation*}
and the relation between the stream function and the temperature is given by
\begin{equation*}
\theta=-(-\Delta)^{1/2}\psi.
\end{equation*}
This system have been considered in frontogenesis, where the dynamics of hot and cold fluids is
studied together with the formation and the evolution of fronts (see \cite{CMT}, \cite{CNS},
\cite{CFR}, \cite{OY}).

From a mathematical point of view, this equation have been presented as a two dimensional model of
the 3-D Euler equation due to their strong analogies (see \cite{CMT}), being the formation of
singularities for a regular initial data an open problem (see \cite{CMT}, \cite{Diego}, \cite{CF}).
Nevertheless the QG equation has global in time weak solutions due to an extra cancellation (see
\cite{Resnick}). A few sparse results are known about weak solutions of the 2-D and 3-D Euler
equation in its primitive-variable form.

An outstanding kind of weak solutions for the QG equation are those in which the temperature takes
two different values in complementary domains, modelling the evolution of a sharp front as follows
\begin{equation}\label{temperature}
\theta(x_1,x_2,t)=\left\{\begin{array}{cl}
                    \theta_1,& \Omega(t)\\
                    \theta_2,& \R^2\smallsetminus\Omega(t).
                 \end{array}\right.
\end{equation}

In this work we study a problem similar to the 2-D vortex patch problem, where the vorticity of the
2-D Euler equation is given by a characteristic function of a domain, and it is considered the
regularity of the free boundary of such domain. For this equation the vorticity satisfies
\begin{align}
\begin{split}\label{2E}
w_t(x,t)&+u(x,t)\cdot\grad w(x,t)=0,\quad x\in\R^2,\\
w(x,0)&=w_0(x),
\end{split}
\end{align}
in a weak sense, and the velocity is given by the Biot-Savart law or analogously
\begin{equation*}
u=\gradp\psi,\quad\mbox{and}\quad w=\Delta\psi.
\end{equation*}
Chemin \cite{Ch} proved global-in-time regularity for the free boundary using paradifferential
calculus. A simpler proof can be found in \cite{bertozzi-Constantin} due to Bertozzi and
Constantin.

We point out that in the QG equation, the velocity is determined from the temperature by singular
integral operators (see \cite{St3}) as follows
\begin{equation}\label{velocityriesz}
u=(-R_2\theta,R_1\theta),
\end{equation}
where $R_1$ and $R_2$ are the Riesz transforms, making the system more singular than \eqref{2E}.

Rodrigo \cite{Rodrigo} proposed the problem of the evolution of a sharp front for the QG equation.
He derived the velocity on the free boundary in the normal direction, and proved local-existence
and uniqueness for a periodic $C^{\infty}$ front, i.e.
\begin{equation*}
\theta(x_1,x_2,t)=\left\{\begin{array}{cl}
                    \theta_1,& \{f(x_1,t)>x_2\}\\
                    \theta_2,& \{f(x_1,t)\leq x_2\},
                 \end{array}\right.
\end{equation*}
with $f(x_1,t)$ periodic, using the Nash-Moser iteration.

In this paper we study a family of contour dynamics equation given by weak solutions of the
following system

\begin{equation}\label{fsystem}
\begin{array}{l}
\theta_t+u\cdot\grad \theta=0,\quad x\in\R^2,\\
\\
u=\gradp \psi,\quad \theta=-(-\Delta)^{1-\al/2}\psi,\quad 0<\al\leq 1,
\end{array}
\end{equation}
where the active scalar $\theta(x,t)$ satisfies \eqref{temperature}. We notice that the case
$\al=0$ is the 2-D vortex patch problem, and $\al=1$ correspond to the sharp front for the QG
equation.

This system was introduced by C\'ordoba, Fontelos, Mancho and Rodrigo in \cite{CFMR}, where they
present a proof of local-existence for a periodic $C^{\infty}$ front, and show evidence of
singularities in finite time. The singular scenario is due to two patches collapse point-wise.

Here we give a proof of local-existence of the system \eqref{fsystem} where the solution satisfies
\eqref{temperature}, with the boundary $\partial\Omega(t)$ given by the curve
\begin{equation*}
\partial\Omega(t)=\{x(\g,t)=(x_1(\g,t),x_2(\g,t))\,:\,\g\in[-\pi,\pi]\},
\end{equation*}
and $x(\g,t)$ belongs to a Sobolev space. In the cases $0<\al<1$ we show uniqueness.

It is well-known (see \cite{Hou} and \cite{Rodrigo}) that in this kind of contour dynamics
equations, the velocity in the tangential direction only moves the particles on the boundary.
Therefore we do not alter the shape of the contour if we change the tangential component of the
velocity; i.e., we are making a change on the parametrization. In the most singular case, $\al=1$
or the QG equation, we need to change the velocity in the tangential direction in order to get
existence in the Sobolev spaces. We take a tangential velocity in such a way that $|\dg x(\g,t)|$
satisfies
$$|\dg x(\g,t)|^2=A(t),$$ and does not depend on $\g$. We would like to cite the work of Hou,
Lowengrub and Shelley \cite{Hou} in which this idea was used to study a contour dynamics problem.

We notice that in order to get a non-singular normal velocity of the curve for $0<\al\leq 1$ (see
\cite{CFMR} and \cite{Rodrigo}), we need a one to one curve, and parameterized in such a way that
$$|\dg x(\g,t)|^2>0.$$ Rigorously, we need that
\begin{equation}\label{cbd}
\frac{|x(\g,t)-x(\g-\e,t)|}{|\e|}>0,\quad \forall\,\g,\e\in[-\pi,\pi],
\end{equation}
therefore we give an initial data satisfying this property, and we prove that this condition is
satisfied locally in time. We point out the importance to take into account the evolution of this
quantity due to the numerical simulations in \cite{CFMR}.

Finally, I wish to thank Antonio C\'ordoba and my thesis advisor Diego C\'ordoba for their strong
influence in this work, their advices and suggestions. The author was partially supported by the
grants PAC-05-005-2 of the JCLM (Spain) and MTM2005-05980 of the MEC (Spain).


\section{The Contour Equation}

In this section we deduce the family of contour equations in term of the free boundary $x(\g,t)$.
We consider the equations given by the system \eqref{QG},
with a velocity satisfying
\begin{equation}\label{velocity}
u(x,t)=\gradp \psi(x,t),
\end{equation}
for the stream function it follows
\begin{equation}\label{stream}
\theta=-(-\Delta)^{1-\al/2}\psi,
\end{equation}
and the active scalar fulfills

\begin{equation}\label{temperature2}
\theta(x_1,x_2,t)=\left\{\begin{array}{cl}
                    \theta_1,& \Omega(t)\\
                    \theta_2,& \R^2\smallsetminus\Omega(t).
                 \end{array}\right.
\end{equation}
 The boundary of $\Omega(t)$ is given by the curve
\begin{equation*}
\partial\Omega(t)=\{x(\g,t)=(x_1(\g,t),x_2(\g,t))\,:\,\g\in[-\pi,\pi]=\T\},
\end{equation*}
with $x(\g,t)$ one to one. Due to the identity \eqref{temperature2}, we find that
\begin{equation*}
\gradp \theta=(\theta_1-\theta_2)\,\dg x(\g,t)\,\delta(x-x(\g,t)),
\end{equation*}
where $\delta$ is the Dirac distribution. Using \eqref{velocity} and \eqref{stream}, we got that
\begin{equation*}
u=-(-\Delta)^{\al/2-1}\gradp \theta.
\end{equation*}
Due to the integral operators $-(-\Delta)^{\al/2-1}$ are Riesz potentials (see \cite{St3}), using
the last to identities we obtain that
\begin{equation}\label{velocityd}
u(x,t)=-\frac{\,\,\Theta_{\al}}{2\pi}\int_{\T}\frac{\dg x(\g-\e,t)}{|x-x(\g-\e,t)|^{\al}}d\e,\\
\end{equation}
for $x\neq x(\g,t)$, and $\Theta_{\al}=(\theta_1-\theta_2)\G(\al/2)/2^{1-\al}\G(2-\alpha/2)$. We
notice that for $\al=1$, if $x\rightarrow x(\g,t)$ the integral in \eqref{velocityd} is divergent.
As we have showed before, we are interested in the normal velocity of the systems. Then we have
that using the identity \eqref{velocityd}, and taking the limit as follows
\begin{equation}\label{limitenormal} u(x,t)\cdot \dg^{\bot}x(\g,t),\quad x\rightarrow x(\g,t),
\end{equation} we obtain
\begin{equation}\label{velocityn}
u(x(\g,t),t)\cdot \dg^{\bot}x(\g,t)=-\frac{\,\,\Theta_{\al}}{2\pi}\int_{\T}\frac{\dg
x(\g-\e,t)\cdot\dg^{\bot}x(\g,t)}{|x(\g,t)-x(\g-\e,t)|^{\al}}d\e.\\
\end{equation}
This identity is well defined for $0<\al\leq 1$ and a one to one curve $x(\g,t)$. Due to the fact
that tangential velocity does not change the shape of the boundary, we fix the contour $\al$-patch
equations as follows
\begin{align}
\begin{split}\label{alpha}
\D x_t(\g,t)&=\frac{\,\,\Theta_\al}{2\pi}\int_{\T}\frac{\dg x(\g,t)-\dg x(\g-\e,t)}{|x(\g,t)-x(\g-\e,t)|^{\al}}d\e,
\quad 0<\al\leq 1,\\
x(\g,0)&=x_0(\g).
\end{split}
\end{align}
Seeing the equation \eqref{velocityd}, we show that the velocity in QG presents a logarithmic
divergence in the tangential direction on the boundary. Nevertheless it belongs to $L^{p}(\R^2)$
for $1<p<\infty$, and to the bounded mean oscillation space (see \cite{St3} for the definition of
the BMO space). In QG the velocity is given by \eqref{velocityriesz}, and writing the temperature
in the following way
$$\theta(x,t)=(\theta_1-\theta_2)X_{\Omega(t)}(x)+\theta_2,$$
we find that
$$u(x,t)=(\theta_1-\theta_2)(-R_2(X_{\Omega(t)}),R_1(X_{\Omega(t)})).$$
Using that $X_{\Omega(t)}\in L^{p}(\R^2)$ for $1\leq p \leq \infty$, we conclude de argument. In
particular the energy of the system is conserved due to
$\|u\|_{L^2}(t)=|\theta_1-\theta_2|\,|\Omega(t)|^{1/2}$, and the area of $\Omega(t)$ is constant in
time.


\section{Weak solutions for the $\al$-system}

In this section we show that if $\theta(x,t)$ is defined by \eqref{temperature2} and the curve
$x(\g,t)$ is convected by the normal velocity \eqref{velocityn}, then $\theta(x,t)$ is a weak
solution of the system \eqref{fsystem} and conversely. We give the definition of weak solution
below.

\begin{defi}
The active scalar $\theta$ is a weak solution of the $\al$-system if for any function $\varphi\in
C_c^{\infty}(\R^2\times(0,T))$, we have
\begin{equation}\label{sd}
\int_0^T\int_{\R^2} \theta(x,t)(\dpt\varphi(x,t)+u(x,t)\cdot\grad\varphi(x,t))dxdt=0,
\end{equation}
where the incompressible velocity $u$ is given by \eqref{velocity}, and the stream function
satisfies \eqref{stream}.
\end{defi}

Then
\begin{prop}
If $\theta(x,t)$ is defined by \eqref{temperature2}, and the curve $x(\g,t)$ satisfies \eqref{cbd}
and \eqref{velocityn}, then $\theta(x,t)$ is a weak solution of the $\al$-system. Furthermore, if
$\theta(x,t)$ is a weak solution of the $\al$-system given by \eqref{temperature2}, and $x(\g,t)$
satisfies \eqref{cbd}, then $x(\g,t)$ verifies \eqref{velocityn}.
\end{prop}

Proof: Let $\theta(x,t)$ be a weak solution of the $\al$-system defined by \eqref{temperature2}.
Integrating by parts we have
\begin{align*}
I&=\int_0^T\int_{\R^2}\theta(x,t)\dpt\varphi(x,t)dxdt=\theta_1\int_0^T\!\!\int_{\Omega(t)}\partial_t\varphi(x,t)
dxdt+\theta_2\int_0^T\!\!\int_{\Omega(t)\smallsetminus\R^2}\!\!\partial_t\varphi(x,t)dxdt\\
&=-(\theta_1-\theta_2)\int_0^T\int_\T\varphi(x(\g,t),t)\, x_t(\g,t)\cdot \dg^{\bot} x(\g,t) d\g dt.
\end{align*}
On the other hand, we obtain
\begin{align*}
J&=\int_0^T\int_{\R^2}\theta\, u\cdot\grad\varphi\,
dxdt=\theta_1\int_0^T\int_{\Omega}u\cdot\grad\varphi\,
dxdt+\theta_2\int_0^T\int_{\R^2\smallsetminus\Omega}u\cdot\grad\varphi\,dxdt.
\end{align*}
Taking $$\Omega_1^\ep(t)=\{x\in\Omega\,:\,dist(x,\Omega(t))\geq \ep\},$$ and
$$\Omega_2^\ep(t)=\{x\in\R^2\smallsetminus\Omega\,:\,dist(x,\R^2\smallsetminus\Omega(t))\geq \ep\},$$
we have that $J^\ep\rightarrow J$ if $\ep\rightarrow 0$, where $J^\ep$ is given by
\begin{align*}
J^\ep&=\theta_1\int_0^T\int_{\Omega_1^\ep(t)}u\cdot\grad\varphi\,
dxdt+\theta_2\int_0^T\int_{\Omega_2^\ep(t)}u\cdot\grad\varphi\,dxdt.
\end{align*}
Integrating by part in $J^{\ep}$, using that the velocity is divergence free, and taking the limit
as in \eqref{limitenormal}, we obtain
\begin{align*}
J&=(\theta_1-\theta_2)\int_0^T\int_{\T}\varphi(x(\g,t),t) u(x(\g,t),t)\cdot \dg^{\bot}x(\g,t)d\g
dt\\
&=-(\theta_1-\theta_2)\frac{\,\,\Theta_{\al}}{2\pi}\int_0^T\int_\T\varphi(x(\g,t),t)\Big(\int_{\T}\frac{\dg
x(\g-\e,t)\cdot\dg^{\bot}x(\g,t)}{|x(\g,t)-x(\g-\e,t)|^{\al}}d\e\Big) d\g dt.\\
\end{align*}
We have that $I+J=0$ using \eqref{sd}, and it follows
\begin{equation*}
\int_0^T\int_\T
f(\g,t)\Big(x_t(\g,t)\cdot\dg^{\bot}x(\g,t)+\frac{\,\,\Theta_{\al}}{2\pi}\int_{\T}\frac{\dg
x(\g-\e,t)\cdot\dg^{\bot}x(\g,t)}{|x(\g,t)-x(\g-\e,t)|^{\al}}d\e\Big)d\g dt=0,\\
\end{equation*}
for $f(\g,t)$ periodic in $\g$. We find that \eqref{velocityn} is satisfied. Following the same
arguments it is easy to check that if $x(\g,t)$ satisfies \eqref{velocityn}, then $\theta$ is a
weak solution given by \eqref{temperature2}.

\section{Local well-posedness for $0<\al<1$}

In this section we prove existence and uniqueness for the contour equation in the cases $0<\al<1$.
We denote the Sobolev spaces by $H^k(\T)$, with norms
$$\|x\|_{H^k}^2=\|x\|_{L^2}^2+\|\dg^k x\|^2_{L^2},$$
and the spaces $C^k(\T)$ with
$$\|x\|_{C^k}=\max_{j\leq k}\|\dg^j x\|_{L^{\infty}}.$$
We need that the curve satisfies
\begin{equation}\label{ci}
\frac{|x(\g,t)-x(\g-\e,t)|}{|\e|}>0,\quad \forall\,\g,\e\in[-\pi,\pi],
\end{equation}
then we define
\begin{equation}\label{df}
F(x)(\g,\e,t)=\frac{|\e|}{|x(\g,t)-x(\g-\e,t)|}\quad \forall\,\g,\e\in[-\pi,\pi],
\end{equation}
with $$F(x)(\g,0,t)=\frac{1}{|\dg x(\g,t)|}.$$ The main theorem in this section is the following

\begin{thm}\label{theorem}
Let $x_0(\g)\in H^{k}(\T)$ for $k\geq 3$ with $F(x_0)(\g,\e)<\infty$. Then there exists a time
$T>0$ so that there is a unique solution to \eqref{alpha} for $0<\al<1$ in $C^1 ([0,T];H^k(\T))$
with $x(\g,0)=x_0(\g)$.
\end{thm}

Proof: We can choose $\Theta_\al=2\pi$ without loss of generality, obtaining the following equation
\begin{align}
\begin{split}\label{equation}
\D x_t(\g,t)&=\int_{\T}\frac{\dg x(\g,t)-\dg x(\g-\e,t)}{|x(\g,t)-x(\g-\e,t)|^{\al}}d\e,\quad 0<\al< 1,\\
x(\g,0)&=x_0(\g).
\end{split}
\end{align}
We present the proof for $k=3$, being analogous for $k>3$, using energy estimates (see
\cite{bertozzi-Majda} for more details). We ignore the time dependence to simplify the notation in
some terms. Considering the quantity

\begin{align}
\begin{split}\label{en2}
\int_{\T}x(\g)\cdot x_t(\g)d\g &=\int_{\T}\int_{\T}x(\g)\cdot\frac{\dg x(\g)-\dg x(\e)}{|x(\g)-x(\e)|^{\al}}d\e d\g\\
&=-\int_{\T}\int_{\T}x(\e)\cdot\frac{\dg x(\g)-\dg x(\e)}{|x(\g)-x(\e)|^{\al}}d\e d\g\\
&=\frac12\int_{\T}\int_{\T}\frac{(x(\g)-x(\e))\cdot(\dg x(\g)-\dg x(\e))}{|x(\g)-x(\e)|^{\al}}d\e d\g\\
&=\frac{1}{2(2-\al)}\int_{\T}\int_{\T}\dg|x(\g)-x(\g-\e)|^{2-\al} d\g d\e\\
&=0,
\end{split}
\end{align}
we obtain
\begin{equation}\label{cn2}
\D\dt\|x\|_{L^2}(t)=0.
\end{equation}
We decompose as follows

\begin{align*}
\int_{\T}\dg^3x(\g)\cdot\dg^3x_t(\g)d\g &= I_1+I_2+I_3+I_4,
\end{align*}
where
\begin{align*}
I_1&=\int_\T\int_\T\dg^3x(\g)\cdot \frac{\dg^4 x(\g)-\dg^4 x(\g-\e)}{|x(\g)-x(\g-\e)|^{\al}}d\e d\g,\\
I_2&=3\!\!\int_\T\!\int_\T\dg^3x(\g)\cdot (\dg^3 x(\g)-\dg^3 x(\g-\e))\dg (|x(\g)-x(\g-\e)|^{-\al})d\e d\g,\\
I_3&=3\!\!\int_\T\!\int_\T\dg^3x(\g)\cdot (\dg^2 x(\g)-\dg^2 x(\g-\e))\dg^2 (|x(\g)-x(\g-\e)|^{-\al})d\e d\g,\\
I_4&=\!\!\int_\T\!\int_\T\dg^3x(\g)\cdot (\dg x(\g)-\dg x(\g-\e))\dg^3 (|x(\g)-x(\g-\e)|^{-\al})d\e d\g.\\
\end{align*}
Operating as in \eqref{en2}, the term $I_1$ becomes

\begin{align*}
I_1&=\frac12\int_\T\int_\T(\dg^3x(\g)-\dg^3x(\g-\e))\cdot \frac{\dg^4 x(\g)-\dg^4 x(\g-\e)}{|x(\g)-x(\g-\e)|^{\al}}d\e d\g\\
&=\frac14\int_\T\int_\T \frac{\dg|\dg^3 x(\g)-\dg^3 x(\g-\e)|^2}{|x(\g)-x(\g-\e)|^{\al}}d\e d\g\\
&=\frac\al4\int_\T\int_\T \frac{|\dg^3 x(\g)-\dg^3 x(\g-\e)|^2(x(\g)- x(\g-\e))\cdot (\dg x(\g)-\dg x(\g-\e))}{|x(\g)-x(\g-\e)|^{\al+2}}d\e d\g.\\
\end{align*}
One finds that
\begin{align*}
I_1&\leq\frac\al4\int_\T\int_\T \frac{|\dg^3 x(\g)-\dg^3 x(\g-\e)|^2|\dg x(\g)-\dg x(\g-\e)|}{|x(\g)-x(\g-\e)|^{\al+1}}d\e d\g,\\
\end{align*}
and due to the inequality $|\dg x(\g)-\dg x(\g-\e)||\e|^{-1}\leq \|x\|_{C^2},$ it follows
\begin{align}
\begin{split}\label{ei1}
I_1&\leq\frac\al4\|x\|_{C^2}\int_\T\int_\T |\e|^{-\al}|F(x)(\g,\e)|^{1+\al}|\dg^3 x(\g)-\dg^3 x(\g-\e)|^2d\e d\g\\
&\leq \frac12\|F(x)\|_{L^\infty}^{1+\al}\|x\|_{C^2}\int_{\T}|\e|^{-\al}\int_\T (|\dg^3 x(\g)|^2+|\dg^3 x(\g-\e)|^2)d\g d\e\\
&\leq \|F(x)\|_{L^\infty}^{1+\al}\|x\|_{C^2}\|\dg^3x\|^2_{L^2}\int_{\T}|\e|^{-\al}d\e\\
&\leq C_\al \|F(x)\|_{L^\infty}^{1+\al}\|x\|_{C^2}\|\dg^3x\|^2_{L^2}.
\end{split}
\end{align} As before, we can obtain $I_2=-6I_1$, and it yields

\begin{align}
\begin{split}\label{ei2}
I_2&\leq C_{\al}\|F(x)\|_{L^\infty}^{1+\al}\|x\|_{C^2}\|\dg^3x\|^2_{L^2}.
\end{split}
\end{align}
In order to estimate the term $I_3$, we consider $I_3=J_1+J_2+J_3$, where
\begin{align*}
J_1&=-3\al\int_\T\int_\T\dg^3x(\g)\cdot (\dg^2 x(\g)-\dg^2 x(\g-\e)) \frac{A(\g,\e)}{|x(\g)-x(\g-\e)|^{\al+2}}d\e d\g,\\
J_2&=-3\al\!\!\int_\T\!\int_\T\dg^3x(\g)\cdot (\dg^2 x(\g)-\dg^2 x(\g-\e))\frac{|\dg x(\g)- \dg x(\g-\e)|^2}{|x(\g)-x(\g-\e)|^{\al+2}}d\e d\g,\\
J_3&=3\al(2+\al)\!\!\int_\T\!\int_\T\dg^3x(\g)\cdot (\dg^2 x(\g)-\dg^2 x(\g-\e)) \frac{(B(\g,\e))^2}{|x(\g)-x(\g-\e)|^{\al+4}}d\e d\g,\\
\end{align*} with
$$A(\g,\e)=(x(\g)- x(\g-\e))\cdot (\dg^2 x(\g)-\dg^2 x(\g-\e)),$$
and
$$B(\g,\e)=(x(\g)- x(\g-\e))\cdot(\dg x(\g)-\dg x(\g-\e)).$$ The identity
\begin{equation}\label{tvm}\dg^2 x(\g)-\dg^2 x(\g-\e)=\e\int_0^1 \dg^3 x(\g+(s-1)\e) ds,\end{equation} yields

\begin{align*}
J_1&\leq 3\int_0^1\!\!\!\int_\T\int_\T |\e|\frac{(|\dg^2 x(\g)|+|\dg^2 x(\g-\e)|)|\dg^3 x(\g)||\dg^3 x(\g+(s-1)\e)|}{|x(\g)-x(\g-\e)|^{\al+1}}d\g d\e ds\\
&\leq 3\|F(x)\|_{L^\infty}^{1+\al}\|x\|_{C^2}\int_0^1\!\!\!\int_\T|\e|^{-\al}\int_\T (|\dg^3 x(\g)|^2+|\dg^3 x(\g+(s-1)\e)|^2)d\g d\e ds\\
&\leq C_\al \|F(x)\|_{L^\infty}^{1+\al}\|x\|_{C^2}\|\dg^3 x\|_{L^2}^2.
\end{align*}
Using \eqref{tvm}, we have for $J_2$
\begin{align*}
J_2&=-3\al\!\!\int_0^1\!\!\!\int_\T\!\int_\T |F(x)(\g,\e)|^{2+\al}\frac{|\dg x(\g)- \dg x(\g-\e)|^2}{\e} \frac{\dg^3 x(\g)\cdot\dg^3 x(\g\!+\!(s-1)\e)}{|\e|^{\al}}d\g d\e ds\\
&\leq 3\|F(x)\|_{L^\infty}^{2+\al}\|x\|_{C^2}^2\int_0^1\!\!\!\int_\T|\e|^{-\al}\int_\T (|\dg^3 x(\g)|^2+|\dg^3 x(\g+(s-1)\e)|^2)d\g d\e ds\\
&\leq C_\al \|F(x)\|_{L^\infty}^{2+\al}\|x\|_{C^2}^2\|\dg^3 x\|_{L^2}^2.
\end{align*}
The term $J_3$ is estimated by
\begin{align*}
J_3&\leq9\!\!\int_0^1\!\!\!\int_\T\!\int_\T |\e|\frac{|\dg x(\g)- \dg x(\g-\e)|^2|\dg^3 x(\g)|
|\dg^3 x(\g\!+\!(s-1)\e)|}{|x(\g)-x(\g-\e)|^{\al+2}}d\g d\e ds\\
&\leq C_\al \|F(x)\|_{L^\infty}^{2+\al}\|x\|_{C^2}^2\|\dg^3 x\|_{L^2}^2.
\end{align*}
We get finally
\begin{equation}\label{ei3}
I_3\leq
C_{\alpha}(\|F(x)\|_{L^\infty}^{1+\al}\|x\|_{C^2}+\|F(x)\|_{L^\infty}^{2+\al}\|x\|_{C^2}^2)\|\dg^3
x\|_{L^2}^2.
\end{equation}
We decompose the term $I_4=J_4+J_5+J_6+J_7+J_8$ as follows

\begin{align*}
J_4&=-\al\int_\T\int_\T\dg^3x(\g)\cdot (\dg x(\g)-\dg x(\g-\e)) \frac{C(\g,\e)}{|x(\g)-x(\g-\e)|^{\al+2}}d\e d\g,\\
J_5&=-3\al\!\!\int_\T\!\int_\T\dg^3x(\g)\cdot (\dg x(\g)-\dg x(\g-\e))\frac{D(\g,\e)}{|x(\g)-x(\g-\e)|^{\al+2}}d\e d\g,\\
J_6&=5\al(\al+2)\!\!\int_\T\!\int_\T\dg^3x(\g)\cdot (\dg x(\g)-\dg x(\g-\e)) \frac{A(\g,\e)B(\g,\e)}{|x(\g)-x(\g-\e)|^{\al+4}}d\e d\g,\\
J_7&=5\al(\al+2)\!\!\int_\T\!\int_\T\dg^3x(\g)\cdot (\dg x(\g)-\dg x(\g-\e)) \frac{B(\g,\e)|\dg x(\g)-\dg x(\g-\e)|^2}{|x(\g)-x(\g-\e)|^{\al+4}}d\e d\g,\\
J_8&=-2\al(\al+2)(\al+4)\!\!\int_\T\!\int_\T\dg^3x(\g)\cdot (\dg x(\g)-\dg x(\g-\e)) \frac{(B(\g,\e))^3}{|x(\g)-x(\g-\e)|^{\al+6}}d\e d\g,\\
\end{align*}
with
$$C(\g,\e)=(x(\g)-x(\g-\e))\cdot (\dg^3 x(\g)-\dg^3 x(\g-\e)),$$
$$D(\g,\e)=(\dg x(\g)-\dg x(\g-\e))\cdot (\dg^2 x(\g)-\dg^2 x(\g-\e)).$$
The most singular term is $J_4$, in such a way that
\begin{align*}
J_4&\leq \|F(x)\|_{L^\infty}^{1+\al}\|x\|_{C^2}\int_{\T}|\e|^{-\al}\int_\T |\dg^3 x(\g)||\dg^3 x(\g)-\dg^3 x(\g-\e)|d\g d\e\\
&\leq C_\al\|F(x)\|_{L^\infty}^{1+\al}\|x\|_{C^2}\|\dg^3 x\|^2_{L^2}.
\end{align*}
For $J_5$, we have
\begin{align*}
J_5&\leq 3\|F(x)\|_{L^\infty}^{2+\al}\|x\|_{C^2}^2\int_{\T}|\e|^{-\al}\int_\T |\dg^3 x(\g)||\dg^2 x(\g)-\dg^2 x(\g-\e)|d\g d\e\\
&\leq C_\al\|F(x)\|_{L^\infty}^{2+\al}\|x\|_{C^2}^2\|\dg^2 x\|_{L^2}\|\dg^3 x\|_{L^2}.
\end{align*}
In a similar way, we obtain
\begin{align*}
J_6&\leq 15\|F(x)\|_{L^\infty}^{2+\al}\|x\|_{C^2}^2\int_{\T}|\e|^{-\al}\int_\T |\dg^3 x(\g)||\dg^2 x(\g)-\dg^2 x(\g-\e)|d\g d\e\\
&\leq C_\al\|F(x)\|_{L^\infty}^{2+\al}\|x\|_{C^2}^2\|\dg^2 x\|_{L^2}\|\dg^3 x\|_{L^2},
\end{align*}
and
\begin{align*}
J_7&\leq 15\|F(x)\|_{L^\infty}^{3+\al}\|x\|_{C^2}^3\int_{\T}|\e|^{-\al}\int_\T |\dg^3 x(\g)||\dg x(\g)-\dg x(\g-\e)|d\g d\e\\
&\leq C_\al\|F(x)\|_{L^\infty}^{3+\al}\|x\|_{C^2}^3\|\dg x\|_{L^2}\|\dg^3 x\|_{L^2}.
\end{align*}
For the term $J_8$, we get
\begin{align*}
J_8&\leq 30\|F(x)\|_{L^\infty}^{3+\al}\|x\|_{C^2}^3\int_{\T}|\e|^{-\al}\int_\T |\dg^3 x(\g)||\dg x(\g)-\dg x(\g-\e)|d\g d\e\\
&\leq C_\al\|F(x)\|_{L^\infty}^{3+\al}\|x\|_{C^2}^3\|\dg x\|_{L^2}\|\dg^3 x\|_{L^2},
\end{align*}
and finally it follows
\begin{equation}\label{ei4}
I_4\leq C_{\alpha}(\|F(x)\|_{L^\infty}^{1+\al}\|x\|_{C^2}+\|F(x)\|_{L^\infty}^{2+\al}\|x\|_{C^2}^2+
\|F(x)\|_{L^\infty}^{3+\al}\|x\|_{C^2}^3)\|x\|_{H^3}^2.
\end{equation}
The inequalities \eqref{ei1}, \eqref{ei2}, \eqref{ei3} and \eqref{ei4} yield

\begin{align*}
\D\dt\|\dg^3x\|^2_{L^2}(t)&\leq
C_{\alpha}\|F(x)\|_{L^\infty}^{3+\al}(t)\|x\|_{C^2}^3(t)\|x\|_{H^3}^2(t).
\end{align*}
Due to the identity $\|x\|_{H^3}^2=\|x\|_{L^2}^2+\|\dg^3 x\|_{L^2}^2$ and \eqref{cn2}, we have
\begin{equation*}
\D\dt\|x\|_{H^3}(t)\leq C_{\alpha}\|F(x)\|_{L^\infty}^{3+\al}(t)\|x\|_{C^2}^3(t)\|x\|_{H^3}(t).
\end{equation*}
Finally, using Sobolev inequalities, we obtain
\begin{equation}\label{esn}
\D\dt\|x\|_{H^3}(t)\leq C_{\alpha} \|F(x)\|_{L^\infty}^{3+\al}(t)\|x\|^4_{H^3}(t).
\end{equation}
Notice that if we use energy methods at this point of the proof (see \cite{bertozzi-Majda} to get
the comprehensive argument), we need to regularize the equation \eqref{equation} as follows

\begin{align}
\begin{split}\label{equationr}
\D x^\ep_t(\g,t)&=\phi_{\varepsilon}*\int_{\T}\frac{\dg (\phi_{\varepsilon}*x^\ep(\g,t)-
\phi_{\varepsilon}*x^\ep(\g-\e,t))}{|x^\ep(\g,t)-x^\ep(\g-\e,t)|^{\al}}d\e,\\
x^\ep(\g,0)&=x_0(\g),
\end{split}
\end{align}
where $\phi_{\varepsilon}$ is a regular approximation to the identity. If the inequality \eqref{ci}
is satisfied initially, due to the properties of the regular approximations to the identity, we get
a Picard system as follows
\begin{align*}
x^\ep_t(\g,t)&=G^{\ep}(x^{\ep}(\g,t)),\\
x^\ep(\g,0)&=x_0(\g),
\end{align*}
where $G^{\ep}$ is Lipschitz. Therefore, for any $\ep>0$, we obtain a time of existence $t_\ep$
where \eqref{ci} is fulfilled. The way to have a time of existence of the system \eqref{equationr}
independent of $\ep$ is to find energy estimates with bounds independent of $\ep$. Next, by taking
$\ep\rightarrow 0$, we get solutions of the original equation. In this particular case, we obtain
\begin{equation*}
\D\dt\|x^\ep\|_{H^3}(t)\leq C_{\alpha}\|F(x^\ep)\|_{L^\infty}^{3+\al}(t)\|x^{\ep}\|^4_{H^3}(t),
\end{equation*}
and if we take $\ep\rightarrow 0$, it is possible that $\|F(x^\ep)\|_{L^\infty}\rightarrow \infty$.
In fact, we have an energy estimate that depends on $\ep$ and then the argument fails. We can not
suppose that if the initial data fulfils \eqref{ci}, there exist a time $t>0$ independent of $\ep$
in which \eqref{ci} is satisfied, because just at this moment of the proof we do not have a
well-posed system when $\ep\rightarrow 0$ (the Lipschitz constant of $G^{\ep}$ goes to infinity
when $\ep\rightarrow 0$).

In order to solve this problem, we consider the evolution of the quantity $\|F(x)\|_{L^{\infty}}$.
Taking $p>2$, it follows
\begin{align*}
\D\dt\|F(x)\|^p_{L^p}(t)&=\dt\int_{\T}\int_\T\big(\frac{|\e|}{|x(\g,t)-x(\g-\e,t)|}\big)^{p}d\g d\e\\
&=-p\int_{\T}\int_\T|\e|^p\frac{(x(\g,t)-x(\g-\e,t))\cdot (x_t(\g,t)-x_t(\g-\e,t))}{|x(\g,t)-x(\g-\e,t)|^{p+2}}d\g d\e\\
&\leq p\int_{\T}\int_\T\big(\frac{|\e|}{|x(\g,t)-x(\g-\e,t)|}\big)^{p+1}\frac{|x_t(\g,t)-x_t(\g-\e,t)|}{|\e|}d\g d\e.\\
\end{align*}
We have
\begin{align*}
x_t(\g)-x_t(\g\!-\!\e)&=\int_{\T}\frac{\dg x(\g)-\dg x(\g-\xi)}{|x(\g)-x(\g-\xi)|^{\al}}d\xi
-\int_{\T}\frac{\dg x(\g-\e)-\dg x(\g-\e-\xi)}{|x(\g-\e)-x(\g-\e-\xi)|^{\al}}d\xi\\
&=\int_{\T}\!(\frac{\dg x(\g)-\dg x(\g-\xi)}{|x(\g)-x(\g-\xi)|^{\al}}-\frac{\dg x(\g)-\dg x(\g-\xi)}{|x(\g-\e)-x(\g-\e-\xi)|^{\al}})d\xi\\
&\quad+\int_{\T}\frac{\dg x(\g)-\dg x(\g-\e)+\dg x(\g-\e-\xi)-\dg x(\g-\xi)}{|x(\g-\e)-x(\g-\e-\xi)|^{\al}}d\xi\\
&=I_5+I_6.
\end{align*}
In order to estimate the term $I_5$, we consider the function $f(a)=a^\al$. For $a,b>0$, we obtain
that
\begin{equation}\label{calculouno}
|a^{\al}-b^{\al}|=\al|\int_0^1(sa+(1-s)b)^{\al-1}(a-b)ds|\leq \al(\min\{a,b\})^{\al-1}|a-b|.
\end{equation} One finds
\begin{align*}
I_5&\leq \int_{\T}\!\frac{|\dg x(\g)\!-\!\dg
x(\g\!-\!\xi)|||x(\g)\!-\!x(\g\!-\!\xi)|^{\al}-|x(\g\!-\!\e)\!-\!x(\g\!-\!\e\!-\!\xi)|^{\al}|}{|x(\g)\!-\!x(\g\!-\!\xi)|^{\al}|
x(\g\!-\!\e)\!-\!x(\g\!-\!\e\!-\!\xi)|^{\al}}d\xi\\
&\leq\|F(x)\|^{2\al}_{L^{\infty}}\|x\|_{C^2}\int_\T|\xi|^{1-\al}\Big|
\Big|\frac{x(\g)\!-\!x(\g\!-\!\xi)}{\xi}\Big|^{\al}-
\Big|\frac{x(\g\!-\!\e)\!-\!x(\g\!-\!\e\!-\!\xi)}{\xi}\Big|^{\al}\Big|d\xi.
\end{align*}
Using \eqref{calculouno}, we get
\begin{align*}
I_5&\leq\al\|F(x)\|^{1+\al}_{L^{\infty}}\|x\|_{C^2}\int_\T|\xi|^{1-\al}\Big|
\Big|\frac{x(\g)\!-\!x(\g\!-\!\xi)}{\xi}\Big|-\Big|\frac{x(\g\!-\!\e)\!-\!x(\g\!-\!\e-\!\xi)}{\xi}\Big|\Big|d\xi\\
&\leq\al\|F(x)\|^{1+\al}_{L^{\infty}}\|x\|_{C^2}\int_\T|\xi|^{-\al}(|x(\g)-x(\g\!-\!\e)|+
|x(\g\!-\!\xi)\!-\!x(\g\!-\!\e\!-\!\xi)|)d\xi\\
&\leq2\al\|F(x)\|^{1+\al}_{L^{\infty}}\|x\|^2_{C^2}|\e|\int_\T|\xi|^{-\al}d\xi\\
&\leq C_\al\|F(x)\|^{1+\al}_{L^{\infty}}\|x\|^2_{C^2}|\e|.
\end{align*}
We obtain for $I_6$ that
\begin{align*}
I_6&\leq \int_{\T}\frac{|\dg x(\g)-\dg x(\g-\e)|+|\dg x(\g-\e-\xi)-\dg x(\g-\xi)|}{|x(\g-\e)-x(\g-\e-\xi)|^{\al}}d\xi\\
&\leq C_\al\|F(x)\|^{\al}_{L^{\infty}}\|x\|_{C^2}|\e|.
\end{align*}
The last two estimates show that
\begin{align*}
\D\dt\|F(x)\|^p_{L^p}(t)&\leq pC_\al\|x\|^2_{C^2}(t)\|F(x)\|^{1+\al}_{L^{\infty}}(t)\int_{\T^2}(F(x)(\g,\e,t))^{p+1}d\g d\e\\
&\leq pC_\al\|x\|^2_{C^2}(t)\|F(x)\|^{2+\al}_{L^{\infty}}(t)\|F(x)\|^p_{L^p}(t),
\end{align*}
and therefore
\begin{align*}
\D\dt\|F(x)\|_{L^p}(t)&\leq C_\al\|x\|^2_{C^2}(t)\|F(x)\|^{2+\al}_{L^{\infty}}(t)\|F(x)\|_{L^p}(t).
\end{align*}
Integrating in time it follows
\begin{align*}
\|F(x)\|_{L^p}(t+h)&\leq \|F(x)\|_{L^p}(t) exp\,\big(
C_{\al}\!\int_t^{t+h}\!\!\!\!\!\!\|x\|^2_{C^2}(s)\|F(x)\|^{2+\al}_{L^{\infty}}(s)ds\big),
\end{align*}
and taking $p\rightarrow \infty$ we obtain
\begin{align*}
\|F(x)\|_{L^\infty}(t+h)&\leq \|F(x)\|_{L^\infty}(t) exp\,\big(
C_{\al}\!\int_t^{t+h}\!\!\!\!\!\!\|x\|^2_{C^2}(s)\|F(x)\|^{2+\al}_{L^{\infty}}(s)ds\big).
\end{align*}
In order to estimate the derivative of the quantity $\|F(x)\|_{L^\infty}(t)$, we use the last
inequality getting
\begin{align*}
\D\dt\|F(x)\|_{L^\infty}(t)&=\lim_{h\rightarrow
0}(\|F(x)\|_{L^\infty}(t+h)-\|F(x)\|_{L^\infty}(t))h^{-1}\\
&\leq \|F(x)\|_{L^\infty}(t)\lim_{h\rightarrow 0}( exp\,\big(
C_{\al}\!\int_t^{t+h}\!\!\!\!\!\!\|x\|^2_{C^2}(s)\|F(x)\|^{2+\al}_{L^{\infty}}(s)ds\big)-1)h^{-1}\\
&\leq C_\al\|x\|^2_{C^2}(t)\|F(x)\|^{3+\al}_{L^{\infty}}(t).
\end{align*}
Applying Sobolev inequalities we conclude that
\begin{align}
\begin{split}\label{enlif}
\D\dt\|F(x)\|_{L^\infty}(t)&\leq C_\al\|x\|^2_{H^3}(t) \|F(x)\|^{3+\al}_{L^{\infty}}(t).
\end{split}
\end{align}
This estimate does not give a global in time bound for $\|F(x)\|_{L^\infty}(t)$ in terms of norms
of $x(\g,t)$. Then, adding the estimate \eqref{enlif} to \eqref{esn}, we have
\begin{equation*}
\D\dt(\|x\|_{H^3}(t)+\|F(x)\|_{L^\infty}(t))\leq
C_\al\|F(x)\|^{3+\al}_{L^{\infty}}(t)\|x\|^4_{H^3}(t),
\end{equation*}
and finally
\begin{equation}
\D\dt(\|x\|_{H^3}(t)+\|F(x)\|_{L^\infty}(t))\leq
C_\al(\|x\|_{H^3}(t)+\|F(x)\|_{L^{\infty}}(t))^{7+\al}.
\end{equation}
Integrating, we get
$$
\D\|x\|_{H^3}(t)+\|F(x)\|_{L^\infty}(t)\leq
\frac{\|x_0\|_{H^3}+\|F(x_0)\|_{L^{\infty}}}{\big(1-tC_\al\big(\|x_0\|_{H^3}
+\|F(x_0)\|_{L^{\infty}}\big)^{6+\al}\big)^{\frac{1}{6+\al}}},
$$ with $C_\alpha$ depending on $\al$. Then, using the regularized problem \eqref{equationr}, the same estimate is obtained
for $x^{\ep}$ instead of $x$. Therefore we get to find a time of existence independent of $\ep$,
and taking $\ep\rightarrow 0$, the existence result follows.

Let $x$ and $y$ be two solutions of the equation \eqref{equation} with $x(\g,0)=y(\g,0)$, and
$z=x-y$. One has that
\begin{align*}
\int_\T z(\g)\cdot z_t(\g)d\g&=\int_\T\int_{\T}z(\g)\cdot(\frac{\dg x(\g)-\dg x(\g-\e)}{|x(\g)-x(\g-\e)|^{\al}}-\frac{\dg x(\g)-\dg x(\g-\e)}{|y(\g)-y(\g-\e)|^{\al}})d\e d\g\\
&\quad+\int_{\T}\int_{\T}\frac{z(\g)\cdot(\dg z(\g)-\dg z(\g-\e))}{|y(\g)-y(\g-\e)|^{\al}}d\e d\g \\
&=I_7+I_8.
\end{align*}
The term $I_7$ is estimated using \eqref{calculouno} by

\begin{align*}
I_7 &\leq \int_\T\int_{\T}\frac{|z(\g)||\dg x(\g)-\dg x(\g-\e)|\big||x(\g)-x(\g-\e)|^{\al}-|y(\g)-y(\g-\e)|^{\al}\big|}{|x(\g)-x(\g-\e)|^{\al}|y(\g)-y(\g-\e)|^{\al}}d\e d\g\\
&\leq
\|F(x)\|_{L^{\infty}}^{\al}\|F(y)\|_{L^{\infty}}^{\al}\|x\|_{C^2}\!\!\int_\T\int_{\T}\!|\e|^{1-\al}|z(\g)|
\Big|\Big|\frac{x(\g)\!-\!x(\g\!-\!\e)}{\e}\Big|^{\al}\!\!\!\!-\!\Big|\frac{y(\g)\!-\!y(\g\!-\!\e)}{\e}\Big|^{\al}\Big|d\e d\g\\
&\leq
\|F(x)\|_{L^{\infty}}\|F(y)\|_{L^{\infty}}\|x\|_{C^2}\!\!\int_\T\int_{\T}\!|\e|^{1-\al}|z(\g)|
\Big|\Big|\frac{x(\g)\!-\!x(\g\!-\!\e)}{\e}\Big|\!-\!\Big|\frac{y(\g)\!-\!y(\g\!-\!\e)}{\e}\Big|\Big|d\e d\g\\
&\leq \|F(x)\|_{L^{\infty}}\|F(y)\|_{L^{\infty}}\|x\|_{C^2}\int_\T\int_{\T}|\e|^{-\al}|z(\g)||z(\g)-z(\g-\e)|d\e d\g\\
&\leq C_\al\|F(x)\|_{L^{\infty}}\|F(y)\|_{L^{\infty}}\|x\|_{C^2}\|z\|_{L^2}^2.
\end{align*}
Integrating by parts in $I_8$ yields

\begin{align*}
I_8&=\frac12\int_{\T}\int_{\T}\frac{(z(\g)-z(\g-\e))\cdot(\dg z(\g)-\dg z(\g-\e))}{|y(\g)-y(\g-\e)|^{\al}}d\e d\g \\
&=\frac14\int_{\T}\int_{\T}\frac{\dg(|z(\g)-z(\g-\e)|^2)}{|y(\g)-y(\g-\e)|^{\al}}d\e d\g \\
&=\frac\al4\int_{\T}\int_{\T}\frac{|z(\g)-z(\g-\e)|^2(y(\g)-y(\g-\e))\cdot(\dg y(\g)-\dg
y(\g-\e))}{|y(\g)
-y(\g-\e)|^{\al+2}}d\e d\g \\
&\leq C_\al \|F(y)\|_{L^{\infty}}^{1+\al}\|y\|_{C^2}\|z\|_{L^2}^2.
\end{align*}
Finally we obtain
\begin{equation*}
\D\dt\|z\|^2_{L^2}(t)\leq C(\al,x,F(x),y,F(y))\|z\|^2_{L^2}(t),
\end{equation*}
and using Gronwall inequality we conclude that $z=0$.

\section{Existence for $\al=1$; the QG sharp front}

In this section we prove existence for the QG sharp front in Sobolev spaces. We give the norm of
the Holder space $C^{k,\frac12}(\T)$ by
$$\|x\|_{C^{k,\frac12}}=\|x\|_{C^k}+\max_{\g,\e\in\T}\frac{|\dg^k x(\g)-\dg^k x(\g-\e)|}{|\e|^{1/2}}.$$
In the case of $\al=1$, we have the following equation
\begin{align}
\begin{split}\label{contourQG}
\D x_t(\g,t)&=\frac{\theta_2-\theta_1}{2\pi}\int_{\T}\frac{\dg x(\g,t)-\dg x(\g-\e,t)}{|x(\g,t)-x(\g-\e,t)|}d\e,\\
x(\g,0)&=x_0(\g).
\end{split}
\end{align}
We take $\theta_2-\theta_1=2\pi$ without lost of generality. This equation loses two derivatives,
therefore the technique applied in the last section does not work. Recalling that we are trying to
solve the QG equation in a weak sense, we can modify the system \eqref{contourQG} in the tangential
direction without changing the shape of the front, as far as the curve satisfies
\begin{align*}
\begin{split}
\D x_t(\g,t)\cdot \dg^{\bot}x(\g,t)&=-\int_{\T}\frac{\dg
x(\g-\e,t)\cdot\dg^{\bot}x(\g,t)}{|x(\g,t)-x(\g-\e,t)|}d\e.
\end{split}
\end{align*}
We showed in section 3 that the temperature $\theta(x,t)$ given by \eqref{temperature2} is a weak
solution of the QG equation. Then we propose to modify the equation \eqref{contourQG} as follows

\begin{align}
\begin{split}\label{QGm}
\D x_t(\g,t)&=\int_{\T}\frac{\dg x(\g,t)-\dg x(\g-\e,t)}{|x(\g,t)-x(\g-\e,t)|}d\e+\lambda(\g,t)\dg x(\g,t),\\
x(\g,0)&=x_0(\g).
\end{split}
\end{align}
The parameter $\lambda(\g,t)$ is to get an extra cancellation in such a way that
\begin{equation}\label{cancelacionextra}
\dg x(\g,t)\cdot \dg^2 x(\g,t)=0.
\end{equation}
Given an initial data satisfying \eqref{ci}, we can reparameterize it obtaining that $|\dg
x(\g,0)|^2=1$, and therefore \eqref{cancelacionextra} is fulfilled at $t=0$. We can not have $|\dg
x(\g,t)|^2=1$ for all time, but
\begin{equation}\label{modulovt}
|\dg x(\g,t)|^2=A(t).
\end{equation}
We have

\begin{align*}
\begin{split}
A'(t)&=2\dg x(\g,t)\cdot \dg x_t(\g,t)\\
&=2\dg x(\g,t)\cdot \dg \Big(\int_{\T}\frac{\dg x(\g,t)-\dg x(\g-\e,t)}{|x(\g,t)-x(\g-\e,t)|}d\e
\Big)+2\dg\lambda(\g,t) A(t),
\end{split}
\end{align*}
and therefore
\begin{align}
\begin{split}\label{dla}
\dg\lambda(\g,t)=\frac{A'(t)}{2A(t)}-\frac{1}{A(t)}\dg x(\g,t)\cdot \dg \Big(\int_{\T}\frac{\dg
x(\g,t)-\dg x(\g-\e,t)}{|x(\g,t)-x(\g-\e,t)|}d\e \Big).
\end{split}
\end{align}
Because $\lambda(\g,t)$ has to be periodic, we obtain
\begin{align}
\begin{split}\label{AppA}
\frac{A'(t)}{2A(t)}=\frac{1}{2\pi A(t)}\int_\T\dg x(\g,t)\cdot \dg \Big(\int_{\T}\frac{\dg
x(\g,t)-\dg x(\g-\e,t)}{|x(\g,t)-x(\g-\e,t)|}d\e \Big) d\g.
\end{split}
\end{align}
Using \eqref{AppA} in \eqref{dla}, and integrating in $\g$, one gets the following formula for
$\la(\g,t)$
\begin{align}
\begin{split}\label{la}
\lambda(\g,t)&=\frac{\g+\pi}{2\pi}\int_\T\frac{\dg x(\g,t)}{|\dg x(\g,t)|^2}\cdot \dg \Big(\int_{\T}
\frac{\dg x(\g,t)-\dg x(\g-\e,t)}{|x(\g,t)-x(\g-\e,t)|}d\e \Big) d\g\\
&\quad-\int_{-\pi}^\g \frac{\dg x(\e,t)}{|\dg x(\e,t)|^2}\cdot \de \Big(\int_{\T}\frac{\dg
x(\e,t)-\dg x(\e-\xi,t)}{|x(\e,t)-x(\e-\xi,t)|}d\xi \Big)d\e,
\end{split}
\end{align}
taking $\la(-\pi,t)=\la(\pi,t)=0$. If we consider solutions of the equation \eqref{QGm} with
$\la(\g,t)$ given by \eqref{la}, it is easy to check that
$$
\dt |\dg x(\g,t)|^2=\la(\g,t)\dg |\dg x(\g,t)|^2+\mu(t)|\dg x(\g,t)|^2,
$$
with
$$ \mu(t)=\frac{1}{\pi}\int_\T\frac{\dg x(\g,t)}{|\dg x(\g,t)|^2}\cdot\dg
\Big(\int_{\T}\frac{\dg x(\g,t)-\dg x(\g-\e,t)}{|x(\g,t)-x(\g-\e,t)|}d\e \Big)d\g.
$$ Solving this linear partial differential equation, if \eqref{cancelacionextra} is satisfied
initially, one finds that the unique solution is given by
\begin{align*}
\begin{split}
|\dg x(\g,t)|^2=|\dg x(\g,0)|^{2}+\frac{1}{\pi}\int_0^t\int_\T\dg x(\g,s)\cdot \dg
\Big(\int_{\T}\frac{\dg x(\g,s)-\dg x(\g-\e,s)}{|x(\g,s)-x(\g-\e,s)|}d\e \Big) d\g ds.
\end{split}
\end{align*}
Therefore we obtain \eqref{modulovt}.

The main theorem of this section is

\begin{thm}\label{theoremQG}
Let $x_0(\g)\in H^{k}(\T)$ for $k\geq 3$ with $F(x_0)(\g,\e)<\infty$. Then there exists a time
$T>0$ so that there is a solution to \eqref{QGm} in $C^1 ([0,T];H^k(\T))$ with $x(\g,0)=x_0(\g)$
and $\lambda(\g,t)$ given by \eqref{la}.
\end{thm}

Proof: Being analogous for $k>3$, we give the proof for $k=3$. We have showed before that
\eqref{modulovt} is satisfied if $x(\g,t)$ is a solution of \eqref{QGm}. Then we can rewrite
$\la(\g,t)$ as follows

\begin{align}
\begin{split}\label{ladt}
\lambda(\g,t)&=\frac{\g+\pi}{2\pi A(t)}\int_\T\dg x(\g,t)\cdot \dg \Big(\int_{\T}\frac{\dg x(\g,t)-\dg x(\g-\e,t)}{|x(\g,t)-x(\g-\e,t)|}d\e \Big) d\g\\
&\quad-\frac{1}{A(t)}\int_{-\pi}^\g \dg x(\e,t)\cdot \de \Big(\int_{\T}\frac{\dg x(\e,t)-\dg
x(\e-\xi,t)}{|x(\e,t)-x(\e-\xi,t)|}d\xi \Big)d\e.
\end{split}
\end{align}
We obtain

\begin{align*}
\begin{split}
\int_{\T}x(\g)\cdot x_t(\g)d\g &=\int_{\T}\int_{\T}x(\g)\cdot\frac{\dg x(\g)-\dg x(\g-\e)}{|x(\g)-x(\g-\e)|}d\e d\g+
\int_{\T} \la(\g) x(\g)\cdot\dg x(\g) d\g\\
&=I_1+I_2,
\end{split}
\end{align*}
One finds that $I_1=0$, since

\begin{align*}
\begin{split}
I_1&=\int_{\T}\int_{\T}x(\g)\cdot\frac{\dg x(\g)-\dg x(\e)}{|x(\g)-x(\e)|}d\e
d\g=-\int_{\T}\int_{\T}
x(\e)\cdot\frac{\dg x(\g)-\dg x(\e)}{|x(\g)-x(\e)|}d\e d\g\\
&=\frac12\int_{\T}\int_{\T}\frac{(x(\g)\!-\!x(\e))\cdot(\dg x(\g)\!-\!\dg
x(\e))}{|x(\g)\!-\!x(\e)|}d\e d\g
=\frac{1}{2}\int_{\T}\int_{\T}\dg|x(\g)-x(\g-\e)| d\g d\e\\
&=0.
\end{split}
\end{align*}
For the term $I_2$, one obtains that $I_2\leq \|\la\|_{L^{\infty}}\|x\|_{L^2}\|\dg x\|_{L^2}$, and
\begin{align*}
\begin{split}
\|\la\|_{L^{\infty}}&\leq\frac{2}{A(t)}\int_{\T} |\dg x(\g)|\Big| \dg \int_{\T}\frac{\dg x(\g)-\dg x(\g-\e)}{|x(\g)-x(\g-\e)|}d\e \Big|d\g\\
&\leq \frac{2}{A(t)}\int_{\T} |\dg x(\g)|\int_{\T}\frac{|\dg^2 x(\g)-\dg^2 x(\g-\e)|}{|x(\g)-x(\g-\e)|}d\e d\g\\
&\quad+\frac{2}{A(t)}\int_{\T} |\dg x(\g)|\int_{\T}\frac{|\dg x(\g)-\dg
x(\g-\e)|^2}{|x(\g)-x(\g-\e)|^2}d\e d\g=J_1+J_2.
\end{split}
\end{align*}
Due to $1/A(t)\leq \|F(x)\|^2_{L^{\infty}}(t)$, we have
\begin{align*}
\begin{split}
J_1&\leq 2\|F(x)\|^3_{L^{\infty}}\int_0^1\!\!\int_{\T}\int_{\T}|\dg^3 x(\g+(s-1)\e)||\dg x(\g)|d\g
d\e ds\leq
2\|F(x)\|^3_{L^{\infty}}\|x\|^2_{H^3},\\
\end{split}
\end{align*}
and
\begin{align*}
\begin{split}
J_2&\leq 2\|F(x)\|^4_{L^{\infty}}\|x\|_{C^1}\int_0^1\!\!\int_{\T}\int_{\T}|\dg^2
x(\g+(s-1)\e)|^2d\g d\e ds\leq
2\|F(x)\|^4_{L^{\infty}}\|x\|^3_{H^3}.\\
\end{split}
\end{align*}
Therefore we obtain that
\begin{equation}\label{cn2QG}
\D\dt\|x\|^2_{L^2}(t)\leq C\|F(x)\|^4_{L^{\infty}}(t)\|x\|^5_{H^3}(t).
\end{equation}
We decompose as follows

\begin{align*}
\int_{\T}\dg^3x(\g)\cdot\dg^3x_t(\g)d\g &=\int_{\T}\dg^3x(\g)\cdot\dg^3\Big(\int_{\T}\frac{\dg
x(\g)-\dg x(\g-\e)}{|x(\g)-x(\g-\e)|}d\e\Big) d\g\\
&\quad+\int_{\T}\dg^3x(\g)\cdot\dg^3(\la(\g)\dg x(\g))d\g\\
&=I_3+I_4.
\end{align*}
We take $I_3=J_3+J_4+J_5+J_6$ where

$$J_3=\int_\T\int_\T\dg^3x(\g)\cdot \frac{\dg^4 x(\g)-\dg^4 x(\g-\e)}{|x(\g)-x(\g-\e)|}d\e d\g,$$
$$J_4=3\int_\T\!\int_\T\dg^3x(\g)\cdot (\dg^3 x(\g)-\dg^3 x(\g-\e))\dg (|x(\g)-x(\g-\e)|^{-1})d\e
d\g,$$
$$J_5=3\int_\T\!\int_\T\dg^3x(\g)\cdot (\dg^2 x(\g)-\dg^2 x(\g-\e))\dg^2 (|x(\g)-x(\g-\e)|^{-1})d\e
d\g,$$
$$J_6=\int_\T\!\int_\T\dg^3x(\g)\cdot (\dg x(\g)-\dg x(\g-\e))\dg^3 (|x(\g)-x(\g-\e)|^{-1})d\e
d\g.$$ The term $J_3$ can be written as

\begin{align*}
J_3&=\frac12\int_\T\int_\T(\dg^3x(\g)-\dg^3x(\g-\e))\cdot \frac{\dg^4 x(\g)-\dg^4 x(\g-\e)}{|x(\g)-x(\g-\e)|}d\e d\g\\
&=\frac14\int_\T\int_\T \frac{\dg|\dg^3 x(\g)-\dg^3 x(\g-\e)|^2}{|x(\g)-x(\g-\e)|}d\e d\g\\
&=\frac14\int_\T\int_\T \frac{|\dg^3 x(\g)-\dg^3 x(\g-\e)|^2(x(\g)- x(\g-\e))\cdot (\dg x(\g)-\dg x(\g-\e))}{|x(\g)-x(\g-\e)|^{3}}d\e d\g.\\
\end{align*}
If we define
$$B(\g,\e)=(x(\g)- x(\g-\e))\cdot(\dg x(\g)-\dg x(\g-\e)),$$ due to \eqref{cancelacionextra}, we obtain that
\begin{align*}
J_3&=\frac14\int_\T\int_\T |F(x)(\g,\e)|^{3}|\dg^3 x(\g)-\dg^3 x(\g-\e)|^2
\frac{B(\g,\e)\e^{-2}-\dg x(\g)\cdot\dg^2x(\g)}{|\e|}d\e d\g.\\
\end{align*}
Using that
$$
\Big|\frac{B(\g,\e)\e^{-2}-\dg x(\g)\cdot\dg^2x(\g)}{\e}\Big|\leq
2\|x\|^2_{C^{2,\frac12}}|\e|^{-1/2},
$$
we find
\begin{align}
\begin{split}\label{ej3QG}
J_3&\leq \|F(x)\|^{3}_{L^{\infty}}\|x\|^2_{C^{2,\frac12}}\int_\T|\e|^{-1/2}\int_\T (|\dg^3 x(\g)|^2+|\dg^3 x(\g-\e)|^2)d\g d\e\\
&\leq C \|F(x)\|^{3}_{L^{\infty}}\|x\|^2_{C^{2,\frac12}}\|\dg^3x\|_{L^2}^2\\
&\leq C \|F(x)\|^{3}_{L^{\infty}}\|x\|_{H^3}^4.
\end{split}
\end{align}
We obtain that $J_4=-6J_3$, and it yields

\begin{align}
\begin{split}\label{ej4QG}
J_4&\leq C\|F(x)\|_{L^\infty}^{3}\|x\|^4_{H^3}.
\end{split}
\end{align}
In order to estimate the term $J_5$, we consider $J_5=K_1+K_2+K_3$, where
$$K_1=-3\int_\T\int_\T\dg^3x(\g)\cdot (\dg^2 x(\g)-\dg^2 x(\g-\e)) \frac{C(\g,\e)}{|x(\g)-x(\g-\e)|^{3}}d\e
d\g,$$
$$K_2=-3\int_\T\int_\T\dg^3x(\g)\cdot (\dg^2 x(\g)-\dg^2 x(\g-\e))\frac{|\dg x(\g)-
\dg x(\g-\e)|^2}{|x(\g)-x(\g-\e)|^{3}}d\e d\g,$$
$$K_3=9\int_\T\int_\T\dg^3x(\g)\cdot (\dg^2 x(\g)-\dg^2 x(\g-\e)) \frac{(B(\g,\e))^2}{|x(\g)-x(\g-\e)|^{5}}d\e
d\g,$$ with
$$C(\g,\e)=(x(\g)- x(\g-\e))\cdot (\dg^2 x(\g)-\dg^2 x(\g-\e)).$$
The inequality
\begin{equation}|\dg^2 x(\g)-\dg^2 x(\g-\e)||\e|^{-1/2}\leq \|x\|_{C^{2,\frac12}},\end{equation} yields

\begin{align*}
K_1&\leq 3\|F(x)\|_{L^{\infty}}^2\|x\|_{C^{2,\frac12}}\int_0^1\int_\T|\e|^{-1/2}\int_\T |\dg^3 x(\g)||\dg^3 x(\g+(s-1)\e)|d\g d\e ds\\
&\leq C\|F(x)\|_{L^\infty}^{2}\| x\|_{H^3}^3.
\end{align*}
As before, we have for $K_2$ that
\begin{align*}
K_2&\leq C\|F(x)\|_{L^\infty}^{3}\|x\|_{C^2}^2\|\dg^3 x\|_{L^2}^2\leq
C\|F(x)\|_{L^\infty}^{3}\|x\|_{H^3}^4.
\end{align*}
The term $K_3$ is estimated by
\begin{align*}
K_3&\leq C\|F(x)\|_{L^\infty}^{3}\|x\|_{C^2}^2\|\dg^3 x\|_{L^2}^2\leq
C\|F(x)\|_{L^\infty}^{3}\|x\|_{H^3}^4.
\end{align*}
We get finally
\begin{equation}\label{ej5QG}
J_5\leq C\|F(x)\|_{L^\infty}^{3}\| x\|_{H^3}^4.
\end{equation}
We decompose the term $J_6=K_4+K_5+K_6+K_7+K_{8}$ as follows

$$K_4=-\int_\T\int_\T\dg^3x(\g)\cdot (\dg x(\g)-\dg x(\g-\e))
\frac{D(\g,\e)}{|x(\g)-x(\g-\e)|^{3}}d\e d\g,$$
$$K_5=-3\int_\T\int_\T\dg^3x(\g)\cdot (\dg x(\g)-\dg x(\g-\e))\frac{E(\g,\e)}{|x(\g)-x(\g-\e)|^{3}}d\e
d\g,$$
$$K_6=15\int_\T\int_\T\dg^3x(\g)\cdot (\dg x(\g)-\dg x(\g-\e)) \frac{B(\g,\e)C(\g,\e)}{|x(\g)-x(\g-\e)|^{5}}d\e
d\g,$$
$$K_7=15\int_\T\int_\T\dg^3x(\g)\cdot (\dg x(\g)-\dg x(\g-\e)) \frac{B(\g,\e)|\dg x(\g)-\dg
x(\g-\e)|^2}{|x(\g)-x(\g-\e)|^{5}}d\e d\g,$$
$$K_{8}=-30\int_\T\int_\T\dg^3x(\g)\cdot (\dg x(\g)-\dg x(\g-\e)) \frac{(B(\g,\e))^3}{|x(\g)-x(\g-\e)|^{7}}d\e
d\g,$$ with
$$D(\g,\e)=(x(\g)-x(\g-\e))\cdot (\dg^3 x(\g)-\dg^3 x(\g-\e)),$$
$$E(\g,\e)=(\dg x(\g)-\dg x(\g-\e))\cdot (\dg^2 x(\g)-\dg^2 x(\g-\e)).$$
We obtain
\begin{align*}
K_5&\leq 3\|F(x)\|_{L^\infty}^{3}\|x\|^2_{C^2}\|\dg^3 x\|^2_{L^2}\leq
3\|F(x)\|_{L^\infty}^{3}\|x\|^4_{H^3},
\end{align*}
\begin{align*}
K_6&\leq 15\|F(x)\|_{L^\infty}^{3}\|x\|^2_{C^2}\|\dg^3 x\|^2_{L^2}\leq
15\|F(x)\|_{L^\infty}^{3}\|x\|^4_{H^3},
\end{align*}
\begin{align*}
K_7&\leq 15\|F(x)\|_{L^\infty}^{4}\|x\|^3_{C^2}\|\dg^3 x\|_{L^2}\|\dg^2 x\|_{L^2}\leq
15\|F(x)\|_{L^\infty}^{4}\|x\|^5_{H^3},
\end{align*}
and
\begin{align*}
K_{8}&\leq 30\|F(x)\|_{L^\infty}^{4}\|x\|^3_{C^2}\|\dg^3 x\|_{L^2}\|\dg^2 x\|_{L^2}\leq
30\|F(x)\|_{L^\infty}^{4}\|x\|^5_{H^3}.
\end{align*}
For the most singular term, we have
\begin{align*}
K_4&=\int_\T\int_\T\dg^3x(\g)\cdot (\dg x(\g)-\dg x(\g-\e)) \frac{\e\,\dg x(\g)\cdot (\dg^3 x(\g)-\dg^3 x(\g-\e))-D(\g,\e)}{|x(\g)-x(\g-\e)|^{3}}d\e d\g\\
&\quad -\int_\T\int_\T\dg^3x(\g)\cdot (\dg x(\g)-\dg x(\g-\e)) \frac{\e\,\dg x(\g)\cdot (\dg^3 x(\g)-\dg^3 x(\g-\e))}{|x(\g)-x(\g-\e)|^{3}}d\e d\g\\
&=L_1+L_2.
\end{align*}
One finds that
\begin{align*}
L_{1}&\leq \|F(x)\|_{L^\infty}^{3}\|x\|^2_{C^2}\int_\T\int_\T |\dg^3 x(\g)||\dg^3 x(\g)-\dg^3
x(\g-\e)|d\g d\e\leq C\|F(x)\|_{L^\infty}^{3}\|x\|^4_{H^3}.
\end{align*}
The term $L_2$ is decomposed, and it yields

\begin{align*}
L_2&=\int_\T\int_\T\dg^3x(\g)\cdot (\dg x(\g)-\dg x(\g-\e)) \frac{\e\,(\dg x(\g)-\dg x(\g-\e))\cdot\dg^3
x(\g-\e)}{|x(\g)-x(\g-\e)|^{3}}d\e d\g\\
&\quad-\!\int_\T\!\int_\T\dg^3x(\g)\!\cdot\! (\dg x(\g)\!-\!\dg x(\g\!-\!\e))\,\e\,\frac{\dg
x(\g)\cdot\dg^3 x(\g)\!-\!
\dg x(\g\!-\!\e)\cdot\dg^3 x(\g\!-\!\e)}{|x(\g)-x(\g-\e)|^{3}}d\e d\g\\
&=M_1+M_2.
\end{align*}
We estimate the term $M_1$ as follows

\begin{align*}
M_1&\leq \|F(x)\|_{L^\infty}^{3}\|x\|^2_{C^2}\int_\T\int_\T |\dg^3 x(\g)||\dg^3 x(\g-\e)|d\g
d\e\leq \|F(x)\|_{L^\infty}^{3}\|x\|^4_{H^3}.
\end{align*}
Taking the derivative in \eqref{cancelacionextra}, we find that $\dg x(\g)\cdot\dg^3 x(\g)=-|\dg^2
x(\g)|^2$, and we rewrite

\begin{align*}
M_2&=\int_\T\int_\T\dg^3x(\g)\cdot (\dg x(\g)-\dg x(\g-\e))\,\e\,\frac{|\dg^2 x(\g)|^2- |\dg^2
x(\g-\e)|^2}{|x(\g)-x(\g-\e)|^{3}}d\e d\g.
\end{align*}
The inequality

\begin{equation}\label{ip} ||\dg^2 x(\g)|^2-|\dg^2 x(\g-\e)|^2|\leq
2\|x\|_{C^2}|\e|\int_0^1|\dg^3 x(\g+(s-1)\e)| ds,
\end{equation} yields

\begin{align*}
M_2&\leq 2\|F(x)\|_{L^\infty}^{3}\|x\|^2_{C^2}\int_0^1\!\int_\T\int_\T |\dg^3 x(\g)||\dg^3
x(\g+(s-1)\e)|d\g d\e ds\leq C\|F(x)\|_{L^\infty}^{3}\|x\|^4_{H^3}.
\end{align*}
We recall that $K_4=L_1+L_2=L_1+M_1+M_2\leq C\|F(x)\|_{L^\infty}^{3}\|x\|^4_{H^3},$ and finally it
follows

\begin{equation}\label{ej6QG}
J_6\leq C\|F(x)\|_{L^\infty}^{4}\|x\|_{H^3}^5.
\end{equation}
Due to \eqref{ej3QG}, \eqref{ej4QG}, \eqref{ej5QG} and \eqref{ej6QG}, we obtain

\begin{equation}\label{ei3QG}
I_3\leq C\|F(x)\|_{L^\infty}^{4}\|x\|_{H^3}^5.
\end{equation}
We take $I_4=J_7+J_8+J_9+J_{10}$, where

\begin{align*}
J_7=\int_\T\la(\g)\,\dg^3x(\g)\cdot \dg^4x(\g)d\g,\quad\,\,\,&\qquad
J_8=3\int_\T\dg \la(\g)\, |\dg^3 x(\g)|^2 d\g,\\
J_9=3\int_\T\dg^2\la(\g)\, \dg^3x(\g)\cdot \dg^2 x(\g)d\g, &\qquad
J_{10}=\int_\T\dg^3\la(\g)\,\dg^3x(\g)\cdot \dg x(\g)d\g.
\end{align*}
We integrate by parts in the term $J_7$, and we get

\begin{align*}
J_7&=-\frac12\int_{\T}\dg\la(\g)|\dg^3 x(\g)|^2 d\g\leq \frac12 \|\dg\la\|_{L^\infty}\|\dg^3
x\|_{L^2}^2.
\end{align*}
Using \eqref{ladt}, we find that
\begin{align}
\begin{split}\label{dladt}
\dg\lambda(\g,t)&=\frac{1}{2\pi A(t)}\int_\T\dg x(\g,t)\cdot \dg \Big(\int_{\T}\frac{\dg x(\g,t)-\dg x(\g-\e,t)}{|x(\g,t)-x(\g-\e,t)|}d\e \Big) d\g\\
&\quad-\frac{1}{A(t)}\dg x(\g,t)\cdot \dg \Big(\int_{\T}\frac{\dg x(\g,t)-\dg x(\g-\e,t)}{|x(\g,t)-x(\g-\e,t)|}d\e\Big)\\
&=K_9+K_{10}.
\end{split}
\end{align}
The term $K_9$ is estimated as $J_{1}$ and $J_{2}$, obtaining
$$
K_{9}\leq \|F(x)\|_{L^\infty}^{4}\|x\|_{H^3}^3.
$$
We have for $K_{10}$ that
\begin{align*}
K_{10}&\leq \frac{\,\,\,\|x\|_{C^2}}{A(t)}\int_\T\Big(\frac{|\dg^2 x(\g,t)-\dg^2 x(\g-\e,t)|}{|x(\g,t)-x(\g-\e,t)|}+
\frac{|\dg x(\g,t)-\dg x(\g-\e,t)|^2}{|x(\g,t)-x(\g-\e,t)|^2}\Big)d\e \\
&\leq 2\|F(x)\|_{L^\infty}^{4}\|x\|^3_{C^{2,\frac12}}\int_{\T}|\e|^{-1/2}d\e\\
&\leq C\|F(x)\|_{L^\infty}^{4}\|x\|^3_{H^3},
\end{align*}
and therefore

\begin{align}\label{ej7QG}
J_7&\leq C\|F(x)\|_{L^\infty}^{4}\|x\|^5_{H^3}.
\end{align}
Due to the identity $J_{8}=-6J_{7}$, one finds that

\begin{align}\label{ej8QG}
J_8&\leq C\|F(x)\|_{L^\infty}^{4}\|x\|^5_{H^3}.
\end{align}
Using that

\begin{align*}
\begin{split}
\dg^2\lambda(\g,t)&=-\frac{1}{A(t)}\dg^2 x(\g,t)\cdot \dg \Big(\int_{\T}\frac{\dg x(\g,t)-\dg x(\g-\e,t)}{|x(\g,t)-x(\g-\e,t)|}d\e\Big)\\
&\quad-\frac{1}{A(t)}\dg x(\g,t)\cdot \dg^2\Big(\int_{\T}\frac{\dg x(\g,t)-\dg x(\g-\e,t)}{|x(\g,t)-x(\g-\e,t)|}d\e\Big),\\
\end{split}
\end{align*}
one gets
\begin{align*}
J_9&=-\frac{1}{A(t)}\int_{\T}\dg^3 x(\g)\cdot\dg^2 x(\g)\,\,\dg^2 x(\g)\cdot \dg  \Big(\int_{\T}\frac{\dg x(\g)-\dg x(\g-\e)}{|x(\g)-x(\g-\e)|}d\e\Big)d\g\\
&\quad-\frac{1}{A(t)}\int_{\T}\dg^3 x(\g)\cdot\dg^2 x(\g)\,\,\dg x(\g)\cdot \dg^2\Big(\int_{\T}\frac{\dg x(\g)-\dg x(\g-\e)}{|x(\g)-x(\g-\e)|}d\e\Big)d\g\\
&=L_3+L_4.
\end{align*}
Therefore
\begin{align*}
L_3&\leq \frac{\,\,\,\|x\|^2_{C^2}}{A(t)}\int_\T\int_\T|\dg^3 x(\g)|\Big(\frac{|\dg^2
x(\g,t)\!-\!\dg^2 x(\g\!-\!\e,t)|}{|x(\g,t)\!
-\!x(\g\!-\!\e,t)|}+\frac{|\dg x(\g,t)\!-\!\dg x(\g\!-\!\e,t)|^2}{|x(\g,t)\!-\!x(\g-\e,t)|^2}\Big)d\e d\g \\
&\leq \|F(x)\|_{L^\infty}^{4}\|x\|^3_{C^{2}}\int_0^1\int_{\T}\int_{\T}|\dg^3 x(\g)|(|\dg^3
x(\g+(t-1)\e)|+
|\dg^2 x(\g+(t-1)\e)|)d\g d\e ds\\
&\leq C\|F(x)\|_{L^\infty}^{4}\|x\|^5_{H^3}.
\end{align*}
Moreover
\begin{align*}
L_4&=-\frac{1}{A(t)}\int_{\T}\int_{\T}\dg^3 x(\g)\cdot\dg^2 x(\g)\,\,\dg x(\g)\cdot \frac{\dg^3
x(\g)-\dg^3 x(\g-\e)}{|
x(\g)-x(\g-\e)|}d\e d\g\\
&+\frac{2}{A(t)}\int_{\T}\int_{\T}\dg^3 x(\g)\cdot\dg^2 x(\g)\,\,\dg x(\g)\cdot \frac{(\dg^2
x(\g)-\dg^2 x(\g-\e))
B(\g,\e)}{|x(\g)-x(\g-\e)|^{3}}d\e d\g\\
&-\frac{1}{A(t)}\int_{\T}\int_{\T}\dg^3 x(\g)\cdot\dg^2 x(\g)\,\,\dg x(\g)\cdot (\dg x(\g)-\dg
x(\g-\e))\dg^2(|x(\g)-
x(\g-\e)|^{-1})d\e d\g\\
&=M_3+M_4+M_5.
\end{align*}
The terms $M_4$ and $M_5$ are estimated as before, and we obtain
\begin{align*}
M_4+M_5\leq C \|F(x)\|_{L^\infty}^{5}\|x\|^6_{H^3}.
\end{align*}
The most singular term is $M_3$, but we find that
\begin{align*}
M_3&=\frac{1}{A(t)}\int_{\T}\int_{\T}\dg^3 x(\g)\cdot\dg^2 x(\g)\,\,\dg^3 x(\g-\e)\cdot \frac{\dg x(\g)-\dg x(\g-\e)}{|x(\g)-x(\g-\e)|}d\e d\g\\
&\quad-\frac{1}{A(t)}\int_{\T}\int_{\T}\dg^3 x(\g)\cdot\dg^2 x(\g)\,\,\frac{\dg^3 x(\g)\cdot\dg x(\g)-\dg^3 x(\g-\e)\cdot\dg x(\g-\e)}{|x(\g)-x(\g-\e)|}d\e d\g\\
&=N_1+N_2.
\end{align*}
We obtain
\begin{align*}
N_1&\leq \|F(x)\|^3_{L^\infty}\|x\|^2_{C^2}\|\dg^3x\|^2_{L^2}\leq
\|F(x)\|^3_{L^\infty}\|x\|^4_{H^3},
\end{align*} and using \eqref{cancelacionextra}
\begin{align*}
N_2&=\frac{1}{A(t)}\int_{\T}\int_{\T}\dg^3 x(\g)\cdot\dg^2 x(\g)\,\,\frac{|\dg^2 x(\g)|^2-|\dg^2 x(\g-\e)|^2}{|x(\g)-x(\g-\e)|}d\e d\g.\\
\end{align*}
Due to \eqref{ip}, we conclude that
\begin{equation*}
N_2\leq2\|F(x)\|^3_{L^\infty}\|x\|^2_{C^2}\|\dg^3x\|^2_{L^2}
\leq2\|F(x)\|^3_{L^\infty}\|x\|^4_{H^3}.
\end{equation*}
We have $J_9=L_3+L_4=L_3+M_3+M_4+M_5=L_3+N_1+N_2+M_4+M_5$, and therefore
\begin{align}\label{ej9QG}
J_9\leq \|F(x)\|^5_{L^\infty}\|x\|^6_{H^3}.
\end{align}
The identity \eqref{cancelacionextra} yields
$$
J_{10}=-\int_\T \dg^3 \la(\g)\, |\dg^2x(\g)|^2 d\g=2\int_\T \dg^2 \la(\g)\,
\dg^3x(\g)\cdot\dg^2x(\g) d\g =\frac23 J_9,
$$ and therefore
\begin{align}\label{ej10QG}
J_{10}\leq \|F(x)\|^5_{L^\infty}\|x\|^6_{H^3}.
\end{align}
Due to the inequalities \eqref{ej7QG}, \eqref{ej8QG}, \eqref{ej9QG}, and \eqref{ej10QG}, we get
\begin{align*}
I_4&\leq C\|F(x)\|_{L^\infty}^{5}\|x\|_{H^3}^6.
\end{align*}
Using \eqref{ei3QG} and the last estimate, we have

\begin{align*}
\D\dt\|\dg^3x\|^2_{L^2}(t)&\leq C\|F(x)\|_{L^\infty}^{5}(t)\|x\|_{H^3}^6(t).
\end{align*}
This inequality and  \eqref{cn2QG} bound the evolution of the Sobolev norms of the curve as follows
\begin{equation}\label{nh3QG}
\D\dt\|x\|_{H^3}(t)\leq C\|F(x)\|_{L^\infty}^{5}(t)\|x\|^5_{H^3}(t).
\end{equation}

We continue the argument considering the evolution of the quantity $\|F(x)\|_{L^{\infty}}(t)$.
Taking $p>2$, it yields
\begin{align*}
\D\dt\|F(x)\|^p_{L^p}(t)&\leq
p\int_{\T}\int_\T\big(\frac{|\e|}{|x(\g,t)-x(\g-\e,t)|}\big)^{p+1}\frac{|x_t(\g,t)-
x_t(\g-\e,t)|}{|\e|}d\g d\e.\\
\end{align*}
We have
\begin{align*}
x_t(\g)-x_t(\g\!-\!\e)&=\int_{\T}\!(\frac{\dg x(\g)-\dg x(\g-\xi)}{|x(\g)-x(\g-\xi)|}-\frac{\dg
x(\g)-\dg
x(\g-\xi)}{|x(\g-\e)-x(\g-\e-\xi)|})d\xi\\
&\quad+\int_{\T}\frac{\dg x(\g)-\dg x(\g-\e)+\dg x(\g-\e-\xi)-\dg x(\g-\xi)}{|x(\g-\e)-x(\g-\e-\xi)|}d\xi\\
&\quad+(\la(\g)-\la(\g-\e))\dg x(\g)+\la(\g-\e)(\dg x(\g)-\dg x(\g-\e))\\
&=I_5+I_6+I_7+I_8.
\end{align*}
The term $I_5$ yields

\begin{align*}
I_5&\leq \int_{\T}\!\frac{|\dg x(\g)\!-\!\dg
x(\g\!-\!\xi)|\,\big||x(\g)\!-\!x(\g\!-\!\xi)|-|x(\g\!-\!\e)\!-\!x(\g\!-\!\e\!-\!\xi)|\,\big|}{|x(\g)\!-\!x(\g\!-\!\xi)||
x(\g\!-\!\e)\!-\!x(\g\!-\!\e\!-\!\xi)|}d\xi\\
&\leq\|F(x)\|^{2}_{L^{\infty}}\|x\|_{C^2}\int_\T|\xi|^{-1}|x(\g)-x(\g-\e)-(x(\g-\xi)-x(\g-\e-\xi))|d\xi\\
&\leq\|F(x)\|^{2}_{L^{\infty}}\|x\|_{C^2}|\e|\int_0^1\int_\T\frac{\big|\dg
x(\g+(s-1)\e)-\dg x(\g+(s-1)\e-\xi)\big|}{|\xi|}d\xi ds\\
&\leq 2\pi\|F(x)\|^{2}_{L^{\infty}}\|x\|^2_{C^2}|\e|.
\end{align*}
For $I_6$ we take
\begin{align*}
I_6&\leq\|F(x)\|_{L^{\infty}}|\e|\int_0^1\int_\T\frac{\big|\dg^2
x(\g+(s-1)\e)-\dg^2 x(\g+(s-1)\e-\xi)\big|}{|\xi|}d\xi ds\\
&\leq\|F(x)\|_{L^{\infty}}\|x\|_{C^{2,\frac12}}|\e|\int_0^1\int_\T|\xi|^{-1/2}d\xi ds\\
&\leq C\|F(x)\|_{L^{\infty}}\|x\|_{C^{2,\frac12}}|\e|
\end{align*}
We have for $I_7$
\begin{align*}
I_7&\leq \frac{\,\,2\|x\|_{C^2}}{A(t)}|\e|\max_{\g}|\dg x(\g)||\dg\Big(\int_{\T}\frac{\dg x(\g)-\dg
x(\g-\e)}{|x(\g)-x(\g-\e)|}\Big)d\e |\\
&\leq 2\|F(x)\|^2_{L^{\infty}}\|x\|^2_{C^2}|\e|\max_{\g}\Big(\int_{\T}\frac{|\dg^2 x(\g)\!-\!\dg^2
x(\g\!-\!\e)|}{|x(\g)\!-\!x(\g\!-\!\e)|}d\e\!+\!\int_{\T}\frac{|\dg x(\g)\!-\!\dg
x(\g\!-\!\e)|^2}{|x(\g)\!-\!x(\g\!-\!\e)|^2}d\e\Big)\\
&\leq 4\|F(x)\|^4_{L^{\infty}}\|x\|^4_{H^3}|\e|.
\end{align*}
Estimating $\|\la\|_{L^{\infty}}$ as before, easily we get

\begin{align*}
I_8\leq \|\la\|_{L^{\infty}}\|x\|_{C^2}|\e|\leq 4\|F(x)\|^4_{L^{\infty}}\|x\|^4_{H^3}|\e|.
\end{align*}
The last four estimates show that

\begin{align*}
\D\dt\|F(x)\|_{L^p}(t)&\leq C\|x\|^4_{H^3}(t)\|F(x)\|^{5}_{L^{\infty}}(t)\|F(x)\|_{L^p}(t),
\end{align*}
by integrating in time and taking $p\rightarrow \infty$, we obtain

\begin{align*}
\|F(x)\|_{L^\infty}(t+h)&\leq \|F(x)\|_{L^\infty}(t) exp\,\big(
C\!\int_t^{t+h}\!\!\!\!\!\!\|x\|^4_{H^3}(s)\|F(x)\|^{5}_{L^{\infty}}(s)ds\big).
\end{align*}
As in the previous section, it follows

\begin{align*}
\begin{split}
\D\dt\|F(x)\|_{L^\infty}(t)&\leq C\|x\|^4_{H^3}(t)\|F(x)\|^{6}_{L^{\infty}}(t).
\end{split}
\end{align*}
Then, due to \eqref{nh3QG} and the above estimate, we find finally that
\begin{equation*}
\D\dt(\|x\|_{H^3}(t)+\|F(x)\|_{L^\infty}(t))\leq C(\|x\|_{H^3}(t)+\|F(x)\|_{L^{\infty}}(t))^{10}.
\end{equation*}
Integrating, we have
$$
\D\|x\|_{H^3}(t)+\|F(x)\|_{L^\infty}(t)\leq
\frac{\|x_0\|_{H^3}+\|F(x_0)\|_{L^{\infty}}}{\big(1-tC\big(\|x_0\|_{H^3}
+\|F(x_0)\|_{L^{\infty}}\big)^{9}\big)^{\frac{1}{9}}},
$$ where $C$ is a constant.

We have used the equality \eqref{cancelacionextra} to obtain the a priori estimates. In order to
get the solution of \eqref{QGm}, we have to choose an appropriate regularized problem preserving
\eqref{cancelacionextra}. We propose the system
\begin{align}
\begin{split}\label{QGred}
\D x^{\ep,\delta}_t(\g,t)&=\phi_{\varepsilon}*\int_{\T}\frac{\dg
(\phi_{\varepsilon}*x^{\ep,\delta}(\g,t)-\phi_{\varepsilon}*x^{\ep,\delta}(\g-\e,t))}{|x^{\ep,\delta}
(\g,t)-x^{\ep,\delta}(\g-\e,t)|+\delta}d\e+\lambda^{\ep,\delta}(\g,t)\dg x^{\ep,\delta}(\g,t),\\
x^{\ep,\delta}(\g,0)&=x_0(\g),
\end{split}
\end{align}
with
\begin{align*}
\begin{split}
\lambda^{\ep,\delta}(\g,t)&=\frac{\g\!+\!\pi}{2\pi}\!\int_\T\frac{\dg x^{\ep,\delta}(\g,t)}{|\dg
x^{\ep,\delta}(\g,t) |^2}\cdot \dg \Big(\phi_{\varepsilon}*\!\!\int_{\T}\frac{\dg
(\phi_{\varepsilon}*x^{\ep,\delta}(\g,t)\!-\!\phi_{\varepsilon}*x^{\ep,\delta}(\g\!-\!\e,t))}{|x^{\ep,\delta}(\g,t)
\!-\!x^{\ep,\delta}(\g\!-\!\e,t)|+\delta}d\e \Big) d\g\\
&\quad-\!\int_{-\pi}^\g \frac{\dg x^{\ep,\delta}(\e,t)}{|\dg x^{\ep,\delta}(\e,t)|^2}\cdot \de
\Big(\phi_{\varepsilon}*\!\!\int_{\T}\frac{\dg(\phi_{\varepsilon}*x^{\ep,\delta}(\e,t)-\phi_{\varepsilon}*x^{\ep,\delta}
(\e-\xi,t))}{|x^{\ep,\delta}(\e,t)-x^{\ep,\delta}(\e-\xi,t)|+\delta}d\xi \Big)d\e.
\end{split}
\end{align*}
We can obtain energy estimates of the system \eqref{QGred} depending on $\ep$ and $\delta$, but
without using \eqref{cancelacionextra}, and therefore we obtain existence of \eqref{QGred}. As long
as the solution exists, we have that
$$\dg x^{\ep,\delta}(\g,t)\cdot \dg^2 x^{\ep,\delta}(\g,t)=0.$$
Using this property of the solution, we obtain energy estimates that depend only on $\delta$, and
taking $\ep\rightarrow 0$ we get a solution of the following equation
\begin{align}
\begin{split}\label{QGrd}
\D x^{\delta}_t(\g,t)&=\int_{\T}\frac{\dg x^{\delta}(\g,t)-\dg x^{\delta}(\g-\e,t))}{|x^{\delta}
(\g,t)-x^{\delta}(\g-\e,t)|+\delta}d\e+\lambda^{\delta}(\g,t)\dg x^{\delta}(\g,t),\\
x^{\delta}(\g,0)&=x_0(\g),
\end{split}
\end{align}
with
\begin{align*}
\begin{split}
\lambda^{\delta}(\g,t)&=\frac{\g\!+\!\pi}{2\pi}\!\int_\T\frac{\dg x^{\delta}(\g,t)}{|\dg
x^{\delta}(\g,t) |^2}\cdot \dg \Big(\int_{\T}\frac{\dg x^{\delta}(\g,t)\!-\!\dg
x^{\delta}(\g\!-\!\e,t)}{|x^{\delta}(\g,t)\!-\!x^{\delta}(\g\!-\!\e,t)|+\delta}d\e \Big) d\g\\
&\quad-\!\int_{-\pi}^\g \frac{\dg x^{\delta}(\e,t)}{|\dg x^{\delta}(\e,t)|^2}\cdot \de
\Big(\int_{\T}\frac{\dg x^{\delta}(\e,t)-\dg x^{\delta}
(\e-\xi,t))}{|x^{\delta}(\e,t)-x^{\delta}(\e-\xi,t)|+\delta}d\xi \Big)d\e.
\end{split}
\end{align*}
Again we have that the solutions of this system satisfy
$$\dg x^{\delta}(\g,t)\cdot \dg^2 x^{\delta}(\g,t)=0,$$ and taking advantage of this, we find energy
estimates independent of $\delta$. If we tend $\delta$ to $0$, we conclude the existence result.


{\small \noindent Francisco Gancedo\\
Instituto de Matem\'aticas y F\'isica Fundamental\\
Consejo Superior de Investigaciones Cient\'ificas\\
Serrano 123, 28006 Madrid, Spain.\\
E-mail address: fgancedo@imaff.cfmac.csic.es}

\end{document}